\documentclass[3p,times]{elsarticle}

\usepackage{ecrc}

\volume{00}

\firstpage{1}

\journalname{Advances in Applied Mathematics}

\runauth{Tian Yanzhao}

\jid{procs}

\jnltitlelogo{\large Advances in Applied Mathematics}

\usepackage{amssymb}

\usepackage[figuresright]{rotating}

\begin{document}

\begin{frontmatter}



\dochead{}

\title{Remarks on the  zero-divisor graph of a commutative ring }


\author{Tian Yanzhao*,Wei Qijiao} \footnote {*Corresponding author.\par
 E-mail addresses:wyjdeyx@sohu.com,qjwei@cuit.edu.cn,\par
 Project supported by the National Natural Science Foundation of China (Grant No. 11471055)}

\address{(School of Mathematics, Chengdu University of Information Technology, Chengdu 610225, China)}

\begin{abstract}
   In 1988, I.Beck showed that the chromatic number of $G(Z_{n})$ is equal to its clique number.In 2004, S.Akbari and A.Mohammadian proved that the edge chromatic number of $G(Z_{n})$ is equal to its maximum degree,in 2008, J.Skowronek-kaziow give formulas calculating the clique number and the maximum degree of $G(Z_{n})$ ,but he have a error about clique number of $G(Z_{n})$ , we consider the zero-divisor graph $G(Z_{n})$ of the ring $Z_n$.we give formulas calculating the clique number of $G(Z_{n})$.We present a constructed method to calculate the clique number. \par
\end{abstract}


\end{frontmatter}


\section{}
\label{}
The concept of zero-divisor graphs of a commutative rings was introduced by I.Beck in 1988\cite{r1},In 1999 \cite{r2},Anderson and Livingston introduced and studied the Zero-divisor graph whose vertices are the non-zero zero-divisors. This graph turns out to best exhibit the properties of the set of zero-divisors of a commutative ring. The zero-divisor graph helps us to study the algebraic properties of rings using graph theoretical tools. We can translate some algebraic properties of a ring to graph theory language and then the geometric properties of graphs help us to explore some interesting results in the algebraic structures of rings. The zero-divisor graph of a commutative ring has been studied extensively by Anderson, Frazier, Lauve, Levy, Livingston and Shapiro, see\cite{r2,r3,r4,r10}. The zero-divisor graph concept has recently been extended to non-commutative rings, see\cite{r5} . A clique in a graph $G$,is a complete subgraph of $G$,the order of the largest clique in a graph $G$ is its clique number \cite{r6},A subgraph $K_m$ with $m$ vertices is called a clique of size $m$ if any two distinct vertices in it are adjacent .The minimum number of colors that can be used to color the edges of $G$ is called the edge chromatic number and is denoted by $\chi_1(G)$.The maximum degree of $G$ is denoted by $\Delta(G).$The vertex chromatic number $\chi(G) $ of a graph $G$,is the minimum $k$ for which $G$ has a $k-vertex$ coloring.The zero-divisor graph of the rings$Z_{n}$ ,denoted by $G(Z_{n})$,is a graph with vertex set in$Z_{n}-\{0\}$,in which two vertices $x$ and $y$ are adjacent if and only if $x\neq y$ and $x\cdot y\equiv 0(mod n)$.\par
In 1988 ,I.Beck showed that the chromatic number of $G(Z_{n})$ is equal to its clique number.In 2004 S.Akbari and A.Mohammadian proved that the edge chromatic number of $G(Z_{n})$ is equal to its maximum degree\cite{r7},in 2008 \cite{r8}J.Skowronek-kaziow give formulas calculating the clique number and the maximum degree of $G(Z_{n})$ ,but he have a error about clique number of $G(Z_{n})$ .\par
For example,$n=420=2^2\cdot 3\cdot 5\cdot 7$,$G(Z_{420})=\{30,42,70,210\}$,the clique number is 4.$n=108=2^2\cdot3^3$,$G(Z_{108})=\{6,18,36,54,72,90\}$.
the clique number is 6.In this paper we give formulas calculating the clique number of $G(Z_{n})$.\par
  let $n=p_{1}^{\alpha_{1}}p_{2}^{\alpha_{2}}\cdots p_{s}^{\alpha_{s}}$ be the prime power factorization of $n$,where $p_{1}<p_{2}\cdots< p_{s}$ are distinct primes and $\alpha_{i}\geq 1$,$s\geq1$.
\section{Results}
In this section,we show  formulas calculating the clique number of $G(Z_{n})$.at the same time ,we give some examples.
\newtheorem{Theorem}{}
\begin{Theorem}
Vizing's Theorem[6,p281]For every nonempty graph $G$ ,then either $\chi_1(G)=\Delta(G)$ or $\chi_1(G)=\Delta(G)+1$.
\end{Theorem}
\begin{Theorem} Theorem 1.The maximal degree in $G(Z_n)$ has the vertex $n/p_1$and the maximum degree is equal to $n/p_1-1$.\end{Theorem}
Proof:This is proved in\cite{r8}.\qed
\begin{Theorem}
Theorem 2. If $n$ is square-free,then the clique number of the graph $G(Z_{n})$ is $s$.if $\alpha_{i}$are even numbers,for all $1\leq i\leq s$,then the clique number is $p_1^{\alpha_1/2}p_2^{\alpha_2/2}\cdots p_s^{\alpha_s/2}-1$,otherwise ,the clique number is
\begin{center} $p_1^{\delta_1/2}\cdots p_r^{\delta_r/2}q_1^{{\beta_1-1}/2}\cdots q_t^{{\beta_t-1}/2}+t-1$,\end{center} where $\delta_i$ is even,$i=1,\ldots, r$ ,$\beta_i$ is odd,$i=1,\ldots ,t$.
\end {Theorem}
Proof:We consider the three case.\par
1.If $n$ is square-free,let $n=p_1p_2\cdots p_s$,where $p_i$ are distinct primes,$1\leq i\leq s$.we consider the set $S=\{\frac{n}{p_s},\frac{n}{p_{s-1}},\ldots ,\frac{n}{p_1}\}$, the product of every pair elements of the set is a multiple of n.i.e. the elements of S is in the vertices set of $G(Z_{n})$,there are no more elements in $G(Z_{n})$,therefore the clique number in this case is equal to $s$.\par
2.If all $ \alpha_i$ are even, $1\leq i\leq s$,then the element $m=p_1^{\alpha_1/2}p_2^{\alpha_2/2}\cdots p_s^{\alpha_s/2}$and element $2m,3m,4m,\ldots,(m-1)m$ form a clique number of $G(Z_{n})$.the element $t$ is the smallest number such that the multiple $(m-1)m$ is possibly the greatest number belonging to $Z_n$ and the clique number in this case is equal to $m-1$.\par
3.If $\alpha_{i}$ are even and odd numbers,let $n= p_1^{\delta_1/2}\cdots p_r^{\delta_r/2}q_1^{\beta_1-1/2}\cdots q_j^{\beta_j-1/2}\cdot q_1\cdots q_j$,$\delta_i$ is even ($i=1,\ldots, r$ ),$\beta_i$ is odd( $i=1,\ldots ,j$) ,$p_1<p_2<\ldots<p_r$,$q_1<q_2<\ldots<q_j$,$\sum_{i=1}^{r}\delta_i+\sum_{i=1}^{j}\beta_i=\sum_{i=1}^{s}\alpha_i$.
We present a constructed method to calculate the clique number.\par
If $j=1$,$n= p_1^{\delta_1/2}\cdots p_r^{\delta_r/2}q_1^{\beta_1-1}\cdot q_1$,let $k=p_1^{\delta_1/2}\cdots p_r^{\delta_r/2}q_1^{(\beta_1-1)/2}$,we consider the set $A=\{k^2,q_1k,2q_1k,3q_1 k,\ldots,q_1(k-1)k\}$,the product of every pair elements of the set $A$ is a multiple of n.i.e. the elements of $A$ are in the vertices set of $G(Z_{n})$,the element $k$ is the smallest number such that the multiple $(k-1)k$ is possibly the greatest number belonging to $Z_n$ and the clique number in this case is equal to $k+1-1=k$.\par
 If $j=2$,$n= p_1^{\delta_1/2}\cdots p_r^{\delta_r/2}q_1^{\beta_1-1} q_2^{\beta_2-1}q_1q_2$.\par
 let $k=p_1^{\delta_1/2}\cdots p_r^{\delta_r/2}q_1^{(\beta_1-1)/2}q_2^{(\beta_2-1)/2}.$we consider the set $B$\par
$B=\{q_1k,q_2k,q_1q_2k,2q_1q_2k,3q_1q_2k,\ldots,q_1q_2(k-1)k\}$,
the product of every pair elements of the set $B$ is a multiple of n.i.e. the elements of $B$ is
in the vertices set of $G(Z_{n})$,the element $k$ is the smallest number such that the multiple $q_1q_2(k-1)k$ is possibly the greatest number belonging to $Z_n$ and the clique number in this case is equal to $k+2-1=k+1$.\par
\begin{center}\ldots \ldots \ldots \ldots \ldots \ldots \ldots \ldots \ldots  \end{center}.\par
If $j=t$,let $c=q_1q_2\cdots q_t$,$n= p_1^{\delta_1/2}\cdots p_r^{\delta_r/2}q_1^{\beta_1-1}\cdots q_t^{\beta_t-1} c$.\par
 let $k=p_1^{\delta_1/2}\cdots p_r^{\delta_r/2}q_1^{(\beta_1-1)/2}\cdots q_t^{(\beta_t-1)/2}.$
we consider the set $C$\par
$C=\{ck/q_1,ck/q_2,\ldots,ck/q_t,ck,2ck,3ck,\ldots,c(k-1)k\}$,
the product of every pair elements of the set $C$ is a multiple of n.i.e. the elements of $C$ is
in the vertices set of $G(Z_{n})$,the element $k$ is the smallest number such that the multiple $c(k-1)k$ is possibly the greatest number belonging to $Z_n$ and the clique number in this case is $k+t-1$.
of course,the number $j\in N$.\par
we conclude that the clique number is equal to \par
\begin{center}$p_1^{\delta_1/2}\cdots p_r^{\delta_r/2}q_1^{(\beta_1-1)/2}\cdots q_t^{(\beta_t-1)/2}+t-1$\end{center} the proof is complete.\qed
\section{Example}
(1)If $n=60=2^2\cdot3\cdot5$,by the theorem, the clique number of $G(Z_{60})$ is equal to 3. the vertices set of $G(Z_{60})$is $G(Z_{60})=\{12,20,30\}$.\par
(2)If $n=2^5\cdot5^3\cdot7^2$,then,by the theorem, the clique number of $G(Z_{196000})$is equal to 141.\par
(3)If $n=3^3\cdot5^2\cdot7^3$,then,by the theorem, the clique number of $G(Z_{231525})$is equal to 106.
\section{Acknowledgments}
 The authors are indebted to the National Natural Science Foundation of China for support.Also the authors thank the referee for her/his valuable comments.




\bibliographystyle{elsarticle-num}



\end{document}